\documentclass[12pt]{article}
\usepackage{amsmath,amssymb}
\usepackage{rotating}
\usepackage{xypic}
\usepackage{amssymb}
\usepackage{amsmath,amscd}
\usepackage{pst-node,pstricks,multido,pst-plot,pst-text,pst-3d}%
\usepackage{graphicx}
\usepackage{amsfonts}
\newtheorem{example}{Example}[section]

\newtheorem{remark}[example]{Remark}

\newtheorem{corollary}[example]{Corollary}

\newtheorem{proposition}[example]{Proposition}

\newtheorem{lemma}[example]{Lemma}

\def\S{{\mathfrak  S}}
\def\cal#1{{\mathfrak #1}}

\def\<{\langle}
\def\>{\rangle}

\def\C{{\mathbb C}}
\def\Z{{\mathbb Z}}
\def\N{{\mathbb N}}
\def\Y{{\mathbb Y}}

\def\A{{\sf A}}

\def\CT{{\rm C.T.}}

\def\ashuff#1#2#3{
\kern 1pt \vrule height#1 \overline{\vrule height#3 width 0pt
\hskip#2} \rule{.3pt}{#1}\overline{\vrule height#3 width 0pt
\hskip#2} \rule{.3pt}{#1} \kern 1pt }

\def\Det{{\rm Det}}

\def\A{{\mathbb A}}\def\X{{\mathbb X}}\def\Y{{\mathbb Y}}

\def\Det{{\rm Det}}

\def\sign{\mbox{sign}}

\title{Hankel hyperdeterminants, rectangular Jack polynomials and even powers of the Vandermonde }
\author{H. Belbachir\footnote{Universit\'e des Sciences et de la
Technologie Houari Boumediene. BP 32 USTHB 16111 Bab-Ezzouar Alger,
Alg\'erie. email: hacenebelbachir@gmail.com}, A.
Boussicault\footnote{Universit\'e de Paris-Est Marne-la-Vallée,
Institut d'\'Electronique et d'Informatique Gaspard-Monge 77454
Marne-la-Vallée Cedex 2. email: Adrien.Boussicault@univ-mlv.fr} and
J.-G. Luque\footnote{Universit\'e de Paris-Est Marne-la-Vallée,
Institut d'\'Electronique et d'Informatique Gaspard-Monge 77454
Marne-la-Vallée Cedex 2. email: Jean-Gabriel.Luque@univ-mlv.fr}}

\begin{document}
\maketitle
 \begin{abstract}
We investigate the link between rectangular Jack polynomials and
Hankel hyperdeterminants. As an application we give an expression of
the even power of the Vandermonde in term of Jack polynomials.
 \end{abstract}
\section{Introduction}
Few after he introduced the modern notation for determinants
\cite{Ca0}, Cayley proposed several extensions to higher dimensional
arrays under the same name {\it hyperdeterminant} \cite{Ca1,Ca2}.
The notion  considered here is apparently the simplest one, defined
for a $k$th order tensor ${\bf M}=(M_{i_1 \cdots i_k})_{1\leq
i_1,\dots,i_k\leq n }$ on an $n$-dimensional space by
$$\Det{\  M}=\frac1{n!}\sum_{\sigma=(\sigma_1,\cdots,\sigma_k)\in\S_n^k}\sign(\sigma)
M^\sigma,$$ where
$\sign(\sigma)=\sign(\sigma_1)\cdots\sign(\sigma_k)$ is the product
of the signs of the permutations,
$M^\sigma=M_{\sigma_1(1)\dots\sigma_k(1)}\cdots
M_{\sigma_1(n)\dots\sigma_k(n)}$ and $\S_n$ denotes the symmetric
group. Note that others hyperdeterminants are found in literature.
For example, those considered by Gelfand, Kapranov and Zelevinsky
\cite{GKZ} are bigger polynomials with geometric properties. Our
hyperdeterminant is a special case of the {\it riciens}
\cite{Rice,Lec} with only alternant indices. For an hypermatrice
with an even number of indices, it generates the space of the
polynomial invariants of lowest degree.
 Easily, one obtains the nullity of $\Det$ when $k$ is odd and
its invariance under the action of the group $SL_n^{\times k}$. Few
references exists on the topics see \cite{Rice,So1,So2,Hau,Bar,Gh}
and before \cite{LT1} the multidimensional analogues of Hankel
determinant do not seem to have been investigated.

In this paper, we discuss about the links between Hankel
hyperdeterminsants and  Jack's symmetric functions indexed by
rectangular partitions. Jack's symmetric functions are a one
parameter (denoted by $\alpha$ in this paper) generalization of
Schur functions.  They were defined by Henry Jack in 1969 in the aim
to interpolate between Schur fonctions ($\alpha=1$) and zonal
polynomials ($\alpha=2$) \cite{J1,J2}.  The story of Jack's
polynomials is closely related to the generalizations of the Selberg
integral \cite{Bar,Kan,Kora,Las,Selb,LT1}. The relation between
Jack's polynomials and hyperdeterminants appeared implicitly in this
context in \cite{LT1}, when one of the author with J.-Y. Thibon gave
an expression of the Kaneko integral \cite{Kan} in terms of Hankel
hyperdeterminant. More recently, Matsumoto computed \cite{Matsu} an
hyperdeterminantal Jacobi-Trudi type formula for rectangular Jack
polynomials.

The paper is organized as follow. In Section \ref{sec1}, after we
recall definitions of Hankel and Toeplitz hyperdeterminants, we
explain that an Hankel hyperdeterminant can be viewed as the {\it
umber} of an even power of the Vandermonde via the substitution
${\int}_\Y: x^n\rightarrow \Lambda^n(\Y)$ where $\Lambda^n(\Y)$
denotes the $n$th elementary symmetric functions on the alphabet
$\Y$. Section \ref{sec2} is devoted to the generalization of the
Matsumoto formula \cite{Matsu} to {\it almost rectangular} Jack
polynomials. In Section \ref{sec3}, we give an equality involving
the substitution ${\int}_\Y$ and skew Jack polynomials. As an
application,   we give in Section \ref{sec5} expressions of even
powers of the Vandermonde determinant in terms of Schur functions
and Jack polynomials.

\section{Hankel and Toeplitz Hyperdeterminants of symmetric functions\label{sec1}}
\subsection{Symmetric functions}
Symmetric functions over an alphabet $\X$ are functions which are
invariant under permutation of the variables. The $\C$-space of the
symmetric functions over $\X$ is an algebra which will be denoted by
$Sym(\X)$.

Let us consider the complete symmetric functions whose generating
series is
$$
\sigma_t(\X):=\sum_i S^i(\X)t^i=\prod_{x\in\X}\frac1{1-xt},
$$
the elementary symmetric functions
$$
\lambda_t(\X):=\sum_i\Lambda^i(\X)t^i=\prod_{x\in\X}({1+xt})=\sigma_{-1}(\X)^{-1},
$$
and power sum symmetric functions
$$
\psi_t(\X):=\sum_i\Psi^i(\X){t^i\over i}=\log(\sigma_t(\X)).
$$

When there is no algebraic relation between the letters of $\X$,
$Sym(\X)$ is a free (associative, commutative) algebra over
complete, elementary or power sum symmetric functions
$$Sym=\C[S^1,S^2,\cdots]=\C[\Lambda^1,\Lambda^2,\cdots]=\C[\Psi^1,\Psi^2,\cdots].$$

 As a consequence,  the algebra
$Sym(\X)$ ($\X$ being infinite or not) is spanned by the set of the
decreasing products of the generators
$$S^\lambda=S^{\lambda_n}\dots S^{\lambda_1},\,
\Lambda^\lambda=\Lambda^{\lambda_n}\dots \Lambda^{\lambda_1},\,
\Psi^\lambda=\Psi^{\lambda_n}\dots \Psi^{\lambda_1},$$ where
$\lambda=(\lambda_1\geq\lambda_{2}\geq\dots\geq\lambda_{n})$ is a
(decreasing) partition.

The algebra $Sym(\X)$ admits also non multiplicative basis. For
example, the monomial functions defined by
$$m_\lambda(\X)=\sum x_{i_1}^{\lambda_1}\cdots x_{i_n}^{\lambda_{n}},$$
 where the sum is over all the distinct
monomials $x_{i_1}^{\lambda_1}\cdots x_{i_n}^{\lambda_{n}}$ with
$x_{i_1},\dots,x_{i_n}\in\X$, and  the Schur functions   defined via
the Jacobi-Trudi formula
\begin{equation}\label{JacTrud1}
S_\lambda(\X)=\det(S^{\lambda_{i}-i+j}(\X)).
\end{equation}
Note that the Schur basis admits other alternative definitions. For
example, it is the only basis such that
\begin{enumerate}
\item It is orthogonal for the scalar product defined on power sums by
\begin{equation}\label{orthoSchur}\langle
\Psi_\lambda,\Psi_\mu\rangle=\delta_{\lambda,\mu}z_\lambda
\end{equation}
 where $\delta_{\mu,\nu}$ is the Kronecker symbol (equal
to $1$ if $\mu=\nu$ and $0$ otherwise) and
$z_\lambda=\prod_ii^{m_i(\lambda)}m_i(\lambda)!$ if $m_i(\lambda)$
denotes the multiplicity of $i$ as a part of $\lambda$.
\item The coefficient of the dominant term in the expansion in the
monomial basis is $1$,
$$
S_\lambda=m_\lambda+\sum_{\mu<\lambda}u_{\lambda\mu}m_\mu.
$$
\end{enumerate}

When $\X=\{x_1,\dots,x_n\}$ is finite, a Schur function has another
determinantal expression
$$
S_\lambda(\X)={\det(x_i^{\lambda_{j}+n-j})\over\Delta(\X)},
$$
where $\Delta(\X)=\prod_{i<j}(x_i-x_j)$ denotes the Vandermonde
determinant.
\subsection{Definitions and General properties}
A Hankel hyperdeterminant is an hyperdeterminant of a tensor whose
entries depend only of the sum of the indices
$M_{i_1,\dots,i_{2k}}=f(i_1+\cdots+i_{2k})$. One of the authors
investigated such polynomials in relation with the Selberg integral
\cite{LT1,LT2}.

Without lost of generality, we will consider the polynomials
\begin{equation}\label{DefHank}
{\cal
H}^k_n(\X)=\Det\left(\Lambda^{i_1+\cdots+i_{2k}}(\X)\right)_{0\leq
i_1,\dots,i_{2k}\leq n-1},
\end{equation}
where $\Lambda_{m}(\X)$ is the $m$th elementary function on the
alphabet $\X$.

Let us consider a shifted version of Hankel hyperdeterminants
\begin{equation}\label{DefShift}
{\cal
H}^k_{v}(\X)=\Det\left(\Lambda^{i_1+\cdots+i_{2k}+v_{i_1}}(\X)\right)_{0\leq
i_1,\dots,i_{2k}\leq n-1}.
\end{equation}
where $v=(v_0,\dots,v_{n-1})\in\Z^n$. Note that (\ref{DefHank})
implies $M_{0,\dots,0}=\Lambda_0(\X)=1$ by convention. But if
$M_{0,\dots,0}\neq 0,1$, this property can be recovered using a
suitable normalization and, if $M_{0\dots0}=0$, by using the shifted
version (\ref{DefShift}) of the Hankel hyperdeterminant.
 As in \cite{Matsu}, one defines To\"eplitz hyperdeterminant by
 giving directly the shifting version
\begin{equation}\label{DefTop}
{\cal
T}^k_{v}(\X)=\Det\left(\Lambda^{i_1+\cdots+i_{k}-(i_{k+1}+\cdots+i_{2k})+v_{i_1}}(\X)\right)_{0\texttt{}\leq
i_1,\dots,i_{2k}\leq n-1}.
\end{equation}
To\"eplitz hyperdeterminants are related to Hankel hyperdeterminants
by the following formulae.
\begin{proposition}
\begin{enumerate}
\item ${\cal H}^k_{v}(\X)=(-1)^{kn(n-1)\over2}{\cal T}^k_{v+(k(n-1))^n}(\X)$
\item  ${\cal T}^k_v(\X)=(-1)^{kn(n-1)\over2}{\cal H}^k_{v+(k(1-n))^n}(\X)$.
\end{enumerate}
\end{proposition}
{\bf Proof} The equalities (1) and (2) are equivalent and are direct
consequences of the definitions (\ref{DefShift}) and (\ref{DefTop}),
$$\begin{array}{rcl}
{\cal
T}_v^{k}(\X)&=&\Det\left(\Lambda^{i_1+\dots+i_k+(n-1-i_{k+1}+\dots+(n-1-i_{2k}-k(n-1))}(\X)\right)\\
&=&(-1)^{kn(n-1)\over2}{\cal H}_{v+(k(1-n))^n}(\X).
\end{array}
$$
$\Box$.

\subsection{The substitution $x^n\rightarrow \Lambda^n(\Y)$}
Let $\X=\{x_1,\dots,x_n\}$ be a finite alphabet and $\Y$ be another
(potentially infinite) alphabet. For simplicity we will denote by
${\int}_\Y$ the substitution
 $$
 {\int}_\Y x^p = \Lambda^p(\Y),
 $$
for each $x\in \X$ and each $p\in\Z$.

The main tool of this paper is the following proposition.

\begin{proposition}\label{D2H} For any integer $k\in\N-\{0\}$, one has
$$
\frac1{n!}{\int}_\Y\Delta(\X)^{2k}={\cal H}^k_{n}(\Y)
$$
where $\Delta(\X)=\prod_{i<j}(x_i-x_j)$.
\end{proposition}
{\bf Proof} It suffices to develop the power of the Vandermonde
determinant
$$\begin{array}{rclcl}
\Delta(\X)^{2k}&=&\det\left(x_i^{j-1}\right)^{2k}
&=&\displaystyle\sum_{\sigma_1,\cdots,\sigma_{2k}\in\S_n}\sign(\sigma_1\cdots\sigma_{2k})
\prod_ix_i^{\sigma_1(i)+\cdots+\sigma_{2k}(i)-2k}.\end{array}
$$
Hence, applying the substitution, one obtains
$$
\begin{array}{rcl}
\displaystyle\frac1{n!}{\int}_\Y\Delta(\X)^{2k}&=&\displaystyle\frac1{n!}\sum_{\sigma_1,\cdots,\sigma_{2k}\in\S_n}\sign(\sigma_1\cdots\sigma_{2k})
\prod_i\Lambda^{\sigma_1(i)+\cdots+\sigma_{2k}(i)-2k}(\Y)\\&=&{\cal
H}^k_{n}(\Y).\end{array}
$$
 $\Box$ \\ \\
More generally, the Jacobi-Trudi formula (\ref{JacTrud1}) implies
the following result.
\begin{proposition}\label{transSchur} One has
$$\frac1{n!}{\int}_\Y
S_{\lambda}(\X)\Delta(\X)^{2k}={\cal H}^k_{{\rm
reverse}_n(\lambda)}(\Y),
$$
where ${\rm reverse}_n(v)=(v_n,\dots,v_1)$ if $v=(v_1,\dots,v_p)$ is
a composition with $p\leq n$ and $v_{p+i}=0$ for $1\leq i\leq n-p$.
\end{proposition}
{\bf Proof} It suffices to remark that
$$
S_{\lambda}(\X)\Delta(\X)^{2k}=\det(x_i^{\lambda_{n-j+1}+j-1})\det\left(x_i^{j-1}\right)^{2k-1},$$
and apply the same computation than in the proof of Proposition
\ref{D2H}. $\Box$
\begin{example}
If $k=1$ then using the second Jacobi-Trudi formula
\begin{equation}\label{JacTrud2}
S_\lambda=\det(\Lambda^{\lambda'_{n-i}-i+j})\end{equation}
 where
$\lambda'$ denotes the conjugate partition of $\lambda$, Proposition
\ref{transSchur} implies
$$\frac1{n!}{\int}_\Y S_{\lambda}(\X)\Delta(\X)^2=(-1)^{n(n-1)\over
2}S_{(\lambda+(n-1)^n)'}(\Y).$$
\end{example}

\section{Jack Polynomials and Hyperdeterminants\label{sec2}}

In this section, we will consider the symmetric functions as a
$\lambda$-ring endowed with the operator $S^i$, and we will use the
definition of addition and multiplication of alphabets in this
context (see {\it e.g.} \cite{Lasc}). Let $\X$ and $\Y$ be two
alphabets, the symmetric functions over the alphabet $\X+\Y$ are
generated by the complete functions $S^i(\X+\Y)$ defined by $
\sigma_t(\X+\Y)=\sigma_t(\X)\sigma_t(\Y)=\sum_iS^i(\X+\Y)t^i. $ If
$\X=\Y$, one has $\sigma_t(2\X):=\sigma_t(\X+\X)=\sigma_t(\X)^2.$
Similarly one defines $\sigma_t(\alpha\X)=\sigma_t(\X)^{\alpha}$. In
particular, the equality $
\sigma_t(-\X)=\prod_x(1-xt)=\lambda_{-t}(\X)$ gives
$S^i(-\X)=(-1)^i\Lambda^i(\X)$. The product of two alphabet $\X$ and
$\Y$ is defined by $ \sigma_t(\X\Y)=\sum
S^i(\X\Y)t^i=\prod_{x\in\X}\prod_{y\in\Y}{1\over 1-xyt}.$ Note that
$\sigma_1(\X\Y)=K(\X,\Y)=\sum_\lambda S_\lambda(\X)S_\lambda(\Y)$ is
the Cauchy Kernel.

\subsection{Jack polynomials}
One considers a one parameter generalization of the scalar product
(\ref{orthoSchur}) defined  by $\langle
\Psi_\lambda,\Psi_\mu\rangle_\alpha=\delta_{\lambda,\mu}z_\lambda
\alpha^{l(\lambda)},$ where $l(\lambda)=n$ denotes the length of the
partition $\lambda=(\lambda_1\geq \dots \geq \lambda_n)$ with
$\lambda_n>0$. The Jack polynomials $P_\lambda^{(\alpha)}$ are the
unique symmetric functions orthogonal for $\langle\, ,\,
\rangle_\alpha$ and such that $
P_\lambda^{(\alpha)}=m_\lambda+\sum_{\mu<\lambda}u_{\lambda\mu}^{(\alpha)}m_\mu.
$ Note that in the case when $\alpha=1$, one recovers the definition
of Schur functions, {\it i.e.} $P_\lambda^{(1)}=S_\lambda$. Let
$(Q_\lambda^{(\alpha)})$ be the dual basis of
$(P_\lambda^{(\alpha)})$. The polynomials $P_\lambda^{(\alpha)}$ and
$Q_\lambda^{(\alpha)}$ are equal up to a scalar factor and the
coefficient of proportionality is computed explicitly in
\cite{Macdo} VI. 10:
\begin{equation}\label{calcscal} b_\lambda^{(\alpha)}:={P_{\lambda}^{(\alpha)}\over  Q_\lambda^{(\alpha)}}=\langle
P_{\lambda}^{(\alpha)}, P_{\lambda}^{(\alpha)}\rangle^{-1}=
\prod_{(i,j)\in\lambda}{\alpha (\lambda_i-j)+\lambda_j'-i+1\over
\alpha (\lambda_i-j+1)+\lambda'_j-i}.
\end{equation}
Let $\X=\{x_1,\cdots,x_n\}$ be a finite alphabet and denote by
$\X^\vee$ the alphabet of the inverse $\{x_1^{-1},\dots,x_n^{-1}\}$.
Let us introduce the second scalar product by
 $$
 \langle f,g\rangle'_{n,\alpha}={1\over
 n!}\CT\{f(\X)g(\X^\vee)\prod_{i\neq
 j}(1-x_ix_j^{-1})^{1\over\alpha}\},
 $$
 see \cite{Macdo} VI. 10.
 The polynomials $P_\lambda^{(\alpha)}$ and
 $Q_\lambda^{(\alpha)}$ are also orthogonal for this scalar product.

For simplicity, we will consider also another normalisation defined
by
$$
R_{\lambda}^{(\alpha),n}(\Y):=\langle
P_{\lambda'}^{(1/\alpha)},Q_{\lambda'}^{(1/\alpha)}\rangle'_{n,\frac
1\alpha}Q^{(\alpha)}_\lambda(\Y).
$$
Note that the polynomial $R^{(\alpha),n}_\lambda$ is not zero only
when $l(\lambda)\leq n$ and in this case the value of the
coefficient $\langle P^{(1/\alpha)}_{\lambda'},
Q^{(1/\alpha)}_{\lambda'}\rangle'_{n,\frac1\alpha}$ is known to be
\begin{equation}\label{calcscalprim}
\begin{array}{l} \langle P^{(1/\alpha)}_{\lambda'},
Q^{(1/\alpha)}_{\lambda'}\rangle'_{n,\frac1\alpha}=\displaystyle
\frac1{n!}\prod_{(i,j)\in\lambda'}{n+\frac1\alpha(j-1)-i+1\over
n+{\frac j\alpha}-i}{\rm C.T.}\left\{\prod_{i\neq
j}\left(1-x_ix_j^{-1}\right)^{\alpha}\right\}
\end{array}
\end{equation}
see \cite{Macdo} VI 10.

\subsection{The operator $\displaystyle{\int}_\Y$ and almost rectangle Jack polynomials}
Suppose that  $\X=\{x_1,\dots,x_n\}$ is a finite alphabet.
² Let $\Y=\{y_1,\cdots\}$ be another (potentially infinite) alphabet
and consider the integral

\begin{equation}\label{depart}
I_{n,k}(\Y)={1\over n!}\CT\{
\Lambda^n(\X^\vee)^{p+k(n-1)}\Lambda^l(\X^\vee)\prod(1+x_iy_j)\Delta(\X)^{2k}\}.
\end{equation}
Note that,
\begin{equation}\label{Ink2}
\begin{array}{l}
I_{n,k}(\Y)= \displaystyle{(-1)^{kn(n-1)\over 2}\over n!}\CT\{
\Lambda^n(\X^\vee)^{p}\Lambda^l(\X^\vee)\prod(1+x_iy_j)\prod_{i\neq
j}(1-{x_i\over x_j})^{k}\}.\end{array}
\end{equation}
But
\begin{equation}\label{LtoP}
\Lambda^n(\X^\vee)^p\Lambda^l(\X^\vee)=P^{(1/k)}_{{(p+1)^lp^{n-l}}}(\X^\vee),
\end{equation}
and
\begin{equation}\label{kernel}
\prod(1+x_iy_j)=\sum_\lambda
Q^{(1/k)}_\lambda(\X)Q^{(k)}_{\lambda'}(\Y).
\end{equation}
Hence, from the orthogonality of $P_\lambda^{(\alpha)}$ and
$Q_\lambda^{(\alpha)}$, equalities (\ref{Ink2}), (\ref{LtoP}) and
(\ref{kernel}) imply
\begin{equation}\label{InkPQ}
I_{n,k}(\Y)={(-1)^{kn(n-1)\over 2}}R^{(k),n}_{n^pl}(\Y).
\end{equation}
On the other hand, one has the equality $${\int}_\Y
x^m=\Lambda_{m}(\Y)=\CT\{x^{-m}\prod_i(1+xy_i)\}. $$ Remarking that
$\Delta(\X^\vee)=(-1)^{n(n-1)\over 2}{\Delta(\X)\over
\Lambda^n(\X)^{n-1}},$ and
$\Lambda^n(\X)^{m}\Lambda^l(\X)=S_{(m^n)+(1^l)}(\X)$ for each
$m\in\Z$, Equality (\ref{depart}) can be written as
\begin{equation}\label{InkSD}
\begin{array}{rclcl}
I_{n,k}(\Y)&=&\displaystyle{1\over n!}{\int}_\Y
S_{((p-k(n-1))^n)+(1^l)}(\X)\Delta(\X)^{2k}&=&
(-1)^{kn(n-1)\over2}{\cal T}^k_{p^{n-l}(p+1)^l}(\X).\end{array}
\end{equation}

One deduces an hyperdeterminantal expression for a Jack polynomial
indexed by the partition $n^pl$.

\begin{proposition}\label{genmatsumoto}
For any positive integers $n, p, l$ and $k$, one has.
$$
R^{(k),n}_{n^pl}={\cal T}_{p^{n-l}(p+1)^l}^k.
$$
\end{proposition}
 The constant term  appearing in (\ref{calcscalprim}) is a special
 case  of the
the Dyson Conjecture \cite{Dyson}. The conjecture of Dyson has been
proved the same year independently by Gunson \cite{Gunson} and
Wilson \cite{Wilson} ( in 1970 I. J. Good \cite{Good} have shown an
elegant elementary proof involving Lagrange interpolation),
$$
{\rm C.T.}\prod_{i\neq
j}\left(1-x_ix_j^{-1}\right)^{a_i}=\left(a_1+\cdots+a_n\atop
a_1,\dots,a_n\right), $$ for $a_1,\dots,a_n\in\N$. Hence, one has
$$
Q^{(k)}_{n^pl}(\Y)={n!}\left(kn\atop k,\cdots,k\right)^{-1}
\kappa(n,p,l;k){\cal T}_{p^{n-l}(p+1)^l}^k(\Y)
$$
where $$\kappa(n,p,l;k)=\prod_{i=1}^n\prod_{j=1}^p{{j}+k(i-1)\over
{j-1}+ki}\prod_{i=1}^l{{p+1}+k(n-i)\over {p}+k(n-i+1)}.$$

 In particular, when $l=0$, one recovers a theorem by Matsumoto.
\begin{corollary} (Matsumoto \cite{Matsu})
$$P^{(k)}_{n^p}(\Y)={ n!}\left(kn\atop
k,\dots,k\right)^{-1}{\cal T}_{p^{n}}^k(\Y).$$
\end{corollary}
{\bf Proof} From equalities (\ref{calcscalprim}) and
(\ref{calcscal}), one has
$$
\langle P^{(1/k)}_{p^n}, Q^{(1/k)}_{p^n}\rangle'_{1/k,n}={ 1\over
n!}\left(kn\atop k,\dots,k\right) \langle P^{(k)}_{n^pl},
P^{(k)}_{n^pl}\rangle_{k}.
$$
Applying Proposition \ref{genmatsumoto}, one finds the result.
 $\Box$

Setting $p=k(n-1)$, one obtains the expression of an Hankel
hyperdeterminant as a Jack polynomials.
\begin{corollary}
$${\cal H}_n^k(\Y)={(-1)^{kn(n-1)\over 2}\over n!}\left(kn\atop
k,\dots,k\right)P^{(k)}_{n^{k(n-1)}}(\Y).$$
\end{corollary}
\subsection{Jack polynomials with parameter $\alpha=\frac1k$}

Let $\Y=\{y_1,y_2,\cdots\}$ be a (potentially infinite alphabet).
 Consider the endomorphism defined on the
power sums symmetric functions $\Psi_p(\Y)$ by
$\omega_\alpha(\Psi_p(\Y)):=\Psi_p(-\alpha\overline\Y)=(-1)^{p-1}\alpha\Psi_p(\Y)$
(see \cite{Macdo} VI 10), where $\overline\Y=\{-y_1,-y_2,\cdots\}$ .
This map is known to satisfy the identities
$$
\omega_\alpha
P_\lambda^{(\alpha)}(\Y)=Q_{\lambda'}^{(\frac1\alpha)}(\Y)
$$
and
$$
\omega_\alpha
\Lambda^n(\Y)=g_{(\frac1\alpha)}^n(\Y):=\Lambda^n(-\alpha\overline\Y).
$$
Applying $\omega_k$ on Proposition \ref{genmatsumoto}, one obtains
the expression of a Jack polynomial with parameter $\alpha=\frac1k$
for an almost rectangular shape $\lambda=(p+1)^lp^{n-l}$ as a
shifted Toeplitz hyperdeterminant whose entries are $$M_{i_1\dots
i_{2k}}=g_{\frac1k}^{i_1+\dots+i_{k}-i_{k+1}-\dots
-i_{2k}+\lambda_{n-i_1+1}}.$$
\begin{proposition} One has
$$
P_{(p+1)^lp^{n-l}}^{(\frac1k)}(\Y)=n!\left(kn\atop
k,\dots,k\right)^{-1}\kappa(n,p,l;k){\cal
T}^{(k)}_{p^{n-l}(p+1)^l}(-k\overline\Y).
$$
\end{proposition}
{\bf Proof} It suffices to apply Proposition \ref{genmatsumoto} with
the alphabet $-\alpha\Y$ to find
$$
Q_{n^pl}^{(k)}(-k\Y)=n!\left(kn\atop
k,\dots,k\right)^{-1}\kappa(n,p,l;k){\cal
T}^{(k)}_{p^{n-l}(p+1)^l}(-k\overline\Y).
$$
The result follows from
$$Q_{n^pl}^{(k)}(-k\overline\Y)=\omega_k Q_{n^pl}^{(k)}(\Y)=P_{(p+1)^lp^{n-l}}^{(\frac1k)}(\Y).$$

$\Box$

\section{Skew Jack polynomials and Hankel hyperdeterminants\label{sec3}}
\subsection{Skew Jack polynomials}
Let us define as in \cite{Macdo} VI 10, the skew $Q$ functions by
$$
\langle Q_{\lambda/\mu}^{(\alpha)}, P_\nu^{(\alpha)}\rangle
:={\langle Q_{\lambda}^{(\alpha)},
P_\mu^{(\alpha)}P_\nu^{(\alpha)}\rangle}.
$$
Straightforwardly, one has
\begin{equation}\label{SQ2Q}
Q^{(\alpha)}_{\lambda/\mu}=\sum_\nu\langle
Q_\lambda^{(\alpha)},P_\nu^{(\alpha)}P_{\mu}^{(\alpha)}\rangle
Q_\nu^{(\alpha)}.
\end{equation}
Classically, the skew Jack polynomials appear when one expands a
Jack polynomial on a sum of alphabet.
\begin{proposition}
Let $\X$ and $\Y$ be two alphabets, one has
 $$Q_{\lambda}^{(\alpha)}(\X+\Y)=\sum_\mu Q^{(\alpha)}_\mu(\X)Q^{(\alpha)}_{\lambda/\mu}(\Y),$$
 or equivalently
$$P_{\lambda}^{(\alpha)}(\X+\Y)=\sum_\mu P^{(\alpha)}_\mu(\X)P^{(\alpha)}_{\lambda/\mu}(\Y).$$
\end{proposition}
{\bf Proof} See \cite{Macdo} VI.7 for a short proof of this
identity. $\Box$\\ \\
 An other
important normalisation is given by
$$J_\lambda^{(\alpha)}=c_\lambda(\alpha)P_\lambda^{(\alpha)}=c'_\lambda(\alpha)Q_\lambda^{(\alpha)},$$
where
$$c_\lambda(\alpha)=\prod_{(i,j)\in\lambda}(\alpha (\lambda_i-j)+\lambda'_j-i+1),$$
and
$$c'_\lambda(\alpha)=\prod_{(i,j)\in\lambda}(\alpha (\lambda_i-j+1)+\lambda'_j-i),$$
if $\lambda'$ denotes the conjugate partition of $\lambda$.

If one defines skew $J$ function by
$$
J^{(\alpha)}_{\lambda/\mu}:=\sum_\nu{\langle
J_\lambda^{(\alpha)},J_\mu^{(\alpha)}J_\nu^{(\alpha)}\rangle_\alpha\over\langle
J_\nu^{(\alpha)},J_\nu^{(\alpha)}\rangle_\alpha} J_\nu^{(\alpha)}
$$
then $J^{(\alpha)}_{\lambda/\mu}$ is again proportional to
$P^{(\alpha)}_{\lambda/\mu}$ and $Q^{(\alpha)}_{\lambda/\mu}$ :
\begin{equation}\label{J2QJ2P}
J^{(\alpha)}_{\lambda/\mu}=c_\lambda(\alpha)c'_\mu(\alpha)P^{(\alpha)}_{\lambda/\mu}=
c'_\lambda(\alpha)c_\mu(\alpha)Q^{(\alpha)}_{\lambda/\mu}.
\end{equation}

\subsection{The operator ${\int}_\Y$ and the skew Jack symmetric functions \label{ssskew}}
Let $\X$, $\Y$ and $\Z$ be three alphabets such that
$\sharp\X=n<\infty$ and $\sharp\Z=m<\infty$.

Consider the polynomial $$ I_{n,k}(\Y,\Z)=\frac1{n!}{\int}_\Y
\prod_ix_i^{-m}\prod_{i,j}(x_i+z_j)\Delta^{2k}(\X).
$$

Remarking that
$$\begin{array}{rcl}
\displaystyle{\int}_\Y
x^{p-m}\prod_i(x+z_i)&=&\displaystyle\sum_{i=0}^m\Lambda^{p+i-m}(\Y)\Lambda^{m-i}(\Z)
\\&=&\displaystyle\Lambda^{p}(\Y+\Z)\\&=&\displaystyle{\int}_{\Y+\Z}x^p,\end{array}
$$
one obtains
\begin{equation}\label{Y+Z}\begin{array}{rcl}
I_{n,k}(\Y,\Z)&=&{\cal
H}^k_{n}(\Y+\Z)\\&=&{(-1)^{kn(n-1)\over2}\over n!}\left(nk\atop
k,\dots,k\right)P^{(k)}_{n^{k(n-1)}}(\Y+\Z).
\end{array}
\end{equation}

Hence, the image of $Q^{\left(1\over
k\right)}_\lambda(\X^\vee)\Delta(\X)^{2k}$ by
$\displaystyle{{\int}}_\Y$ is a Jack polynomial.
\begin{corollary}\label{Q1/XSQ}
One has,
\begin{eqnarray*}
\displaystyle{{\int}}_\Y
 Q^{\left(1\over k\right)}_\lambda(\X^\vee)\Delta(\X)^{2k}=
 \displaystyle{(-1)^{kn(n-1)\over2}}\left(nk\atop
k,\dots,k\right)(b_{n^{k(n-1)}}^{(k)})^{-1}Q_{n^{k(n-1)}/\lambda'}^{(k)}(\Y),
\end{eqnarray*}
 where   $b_{n^{k(n-1)}}^{(k)}={(2(n-1))!(nk)!((n-1)k)!\over k n!(n-1)!((2n-1)k-1)!}.
$

\end{corollary}
{\bf Proof} The equality follows from $$
\prod_ix_i^{-m}\prod_{i,j}(x_i+z_j)=\prod_{i,j}(1+{z_j\over
x_i})=\sum_\lambda Q_{\lambda'}^{(k)}(\Z)Q_{\lambda}^{\left(1\over
k\right)}(\X^\vee).
$$
Indeed, one has
$$
I_{n,k}(\Y,\Z)=\frac1{n!}\sum_\lambda
b_{\lambda'}^{(k)}P_{\lambda'}^{(k)}(\Z){\int}_\Y
Q_\lambda^{\left(1\over k\right)}(\X^\vee) \Delta(\X)^{2k}.
$$
And in the other hand, by (\ref{Y+Z}) one obtains
$$
I_{n,k}(\Y,\Z)={(-1)^{kn(n-1)\over2}\over n!}\left(nk\atop
k,\dots,k\right)\sum_\lambda
P_\lambda^{(k)}(\Z)P^{(k)}_{n^{k(n-1)}/\lambda}(\Y).
$$
Identifying the coefficient of $P_{\lambda'}^{(k)}(\Z)$ in the two
expressions, one finds,
\begin{eqnarray*}
\displaystyle{{\int}}_\Y
 Q^{\left(1\over k\right)}_\lambda(\X^\vee)\Delta(\X)^{2k}=
{(-1)^{kn(n-1)\over2}}\left(nk\atop
k,\dots,k\right)(b_{\lambda}^{(k)})^{-1}P_{n^{k(n-1)}/\lambda'}^{(k)}(\Y),
\end{eqnarray*}
where the value of $b_\lambda^{(k)}:={P_\lambda^{(k)}\over
Q^{(k)}_\lambda}$ is given by equality (\ref{calcscal}). But, from
\ref{J2QJ2P}, $P_{\lambda/\mu}^{(\alpha)}={b_\lambda^{(\alpha)}\over
b_\mu^{(\alpha)}}Q_{\lambda/\mu}^{(\alpha)}$. Hence,
\begin{eqnarray*}
\displaystyle{{\int}}_\Y
 Q^{\left(1\over k\right)}_\lambda(\X^\vee)\Delta(\X)^{2k}=
 \displaystyle{(-1)^{kn(n-1)\over2}}\left(nk\atop
k,\dots,k\right)(b_{n^{k(n-1)}}^{(k)})^{-1}Q_{n^{k(n-1)}/\lambda'}^{(k)}(\Y).
\end{eqnarray*}
The value of $b_{n^{k(n-1)}}^{(k)}$ is obtained from Equality
(\ref{calcscal}) after simplification.
 $\Box$

\section{Even powers of the Vandermonde determinant\label{sec5}}
\subsection{Expansion of the even power of the Vandermonde on the Schur basis}

The expansion of even power of the Vandermonde polynomial on the
Schur functions is an open problem related  the fractional quantum
Hall effect as described by Laughlin's wave function \cite{Lau}. In
particular is  of considerable interest to determine for what
partitions the coefficients of the Schur functions in the expansion
of the square of Vandermonde vanishe \cite{DGIL,Wy1,Wy2,STW,KTW}.
The aim of this subsection is to give an hyperdeterminantal
expression for the coefficient of $S_\lambda(\X)$ in
$\Delta(\X)^{2k}$.

 Let us denote by $\A_{\bf 0}$ the alphabet verifying
$\Lambda^n(\A_{\bf 0})=0$ for each $n\neq 0$ (and by convention
$\Lambda^0(\A_{\bf 0})=1$). The second orthogonality of Jack
polynomials can be written as
$$
\langle f,g\rangle'_{n,\alpha}=\frac1{n!}{\int}_{\bf
\A_0}f(\X)g(\X^\vee)\prod_{i\neq j}(1-x_i^{-1}x_j)^{1\over\alpha}.
$$
In the case when $\alpha=1$, it coincides with the first scalar
product. In particular,
$$
\langle S_\lambda(\X),S_\mu(\X)\rangle'_{n,1}=\delta_{\lambda\mu}.
$$
Hence, the coefficient of $S_\lambda(\X)$ in the expansion of
$\Delta(\X)^{2k}$ is
$$
\langle S_\lambda(\X),\Delta(\X)^{2k}\rangle'_{n,1}=
\frac1{n!}{\int}_{\A_{\bf
0}}S_\lambda(\X)\Delta(\X^{\vee})^{2k}\prod_{i\neq
j}(1-x_i^{-1}x_j).
$$

 One has
$$
\begin{array}{rcl}
\langle S_\lambda(\X),\Delta(\X)^{2k}\rangle'_{n,1}&=&\displaystyle
(-1)^{n(n-1)\over2}\frac1{n!}{\int}_{\A_{\bf
0}}S_\lambda(\X)\Lambda^n(\X)^{(2k+1)(1-n)}\Delta(\X)^{2(k+1)}\\
&=&\displaystyle (-1)^{n(n-1)\over2}\frac1{n!}{\int}_{\A_{\bf
0}}S_{\lambda+((2k+1)(1-n))^n}(\X)\Delta(\X)^{2(k+1)}.
\end{array}
$$
By Proposition \ref{transSchur}, one obtains an hyperdeterminantal
expression for the coefficients of the Schur functions in the
expansion of the even power of the Vandermonde determinant.
\begin{corollary}
The coefficient of $S_\lambda(\X)$ in the expansion of
$\Delta(\X)^{2k}$ is the hyperdeterminant

$$
\langle
S_\lambda(\X),\Delta(\X)^{2k}\rangle'_{n,1}=(-1)^{n(n-1)\over2}{\cal
H}_{{\rm reverse}_n(\lambda)-((2k+1)(n-1))^n}^{k+1}(\A_{\bf 0}).
$$
\end{corollary}
It should be interesting to study the link between the notion of
admissible partitions introduced by Di Francesco and al \cite{DGIL}
and such an hyperdeterminantal expression.

\subsection{Jack polynomials over the alphabet $-\X$ }

In this paragraph, we work with Laurent polynomials in
$\X=\{x_1,\cdots,x_n\}$. The space of symmetric Laurent polynomials
is spanned by the family indexed by decreasing vectors $(\tilde
S_\lambda(\X))_{(\lambda_1\geq\dots\geq \lambda_n)\in\Z^n}$ and
defined by
$$
\tilde S_\lambda(\X)={\det(x_i^{\lambda_j+n-j})\over \Delta(\X)}.
$$
Indeed, each symmetric Laurent polynomial $f$ can be written as
$$
f(\X)=\Lambda^n(\X)^{-m}g(\X)
$$
where $g(\X)$ is a symmetric polynomial in $\X$. As, $g(x)$ is a
linear combination of Schur functions, it follows that $f(\X)$ is a
linear combination of $\tilde S_\lambda$'s.
 Let
$\overline\X=\{-x_1,\dots,-x_n\}$ be the alphabet  of the inverse of
the letters of $\X$.  We consider the operation
${\int}_{-\overline\X}$ ({\it i.e.} the substitution sending each
$x^p$ for $x\in\X$ and $p\in\Z$ to the complete symmetric function
$S^n(\X)$).

 Consider the alternant
$$\begin{array}{rcl}
    a_\lambda(\X)&:=&\displaystyle\sum_\sigma\epsilon(\sigma)x^{\sigma\lambda}\\
    &=& \det(x_i^{\lambda_j})\\
    &=& \tilde S_{\lambda-\delta}(\X)\Delta(\X),\end{array}$$
    where $\delta=(n-1,n-2,\dots, 1,0)$.
From this definition, one obtains that the operator
$\frac1{n!}{\int}_{-\overline \X}$ sends the product of $2k$
alternants $a_\lambda,\,a_\mu,\dots,\,a_\rho$ is an hyperdeterminant

\begin{equation}\label{altodet}
\begin{array}{l}
\displaystyle {1\over n!}{\int}_{-\overline\X}
a_{\lambda}(\X)a_{\mu}(\X)\cdots
a_{\rho}(\X)=\\
\ \ \ \ \ \ \ \ \ \ \ \ \ \ \ \ \ \ \ \ \ \ \ \displaystyle
\Det(S^{\lambda_{i_1}+\mu_{i_2}+\cdots+\rho_{i_{2ki}}}(\X))_{1\leq
i_1,\dots,i_{2k}\leq n}.\end{array}\end{equation}

Consider the linear operator $\Omega^+$ defined by

\begin{equation}\label{Omegaplus}
\Omega^+ \tilde S_\lambda(\X):=
S_\lambda(\X):=\det(S^{\lambda_i+i-j}(\X)).
\end{equation}
In particular, the operator $\Omega^+$ lets invariant the symmetric
polynomials.

Furthermore, it admits an expression involving  ${\int}_{-\overline
\X}$.

\begin{lemma}\label{lemma_alt}
One has
$$
\Omega^+=\frac1{n!}{\int}_{-\overline X}a_\delta(\X)a_{-\delta}(\X).
$$
\end{lemma}
{\bf Proof } It suffices to show that
$$
\frac1{n!}{\int}_{-\overline X}a_\delta(\X)a_{-\delta}(\X)\tilde
S_\lambda(\X)=S_\lambda(\X).
$$
But
$$
a_\delta(\X)a_{-\delta}(\X)\tilde
S_\lambda(\X)=a_{\lambda+\delta}(\X)a_{-\delta}(\X),
$$
and by (\ref{altodet}), one obtains the result.
$\Box$.\\ \\
\begin{proposition}\label{Pat-X}
\ \\    Let $\X=\{x_1,\cdots,x_n\}$ be a finite alphabet and $0\leq
l\leq p\in\N$.
\begin{equation}\label{SVand}
R^{(k),n}_{n^{p+(k-1)(n-1)}l}(-\X)=(-1)^{{(k-1)n(n-1)\over
2}+np+l}S_{(p+1)^lp^{n-l}}(\X)\Delta(\X)^{2(k-1)}.
\end{equation}
\end{proposition}
{\bf Proof} By Lemma \ref{lemma_alt}, one has
$$
\begin{array}{rcl}
S_\lambda(\X)\Delta(\X)^{2(k-1)}&=&\Omega^+S_\lambda(\X)\Delta(\X)^{(2k-1)}\\
&=& \displaystyle\frac1{n!}{\int}_{-\overline\X}
a_{\lambda+\delta}(\X)a_\delta(\X)^{2(k-1)}(\X)a_{-\delta}(\X).
\end{array}
$$
By Equality (\ref{altodet}), one obtains
$$
\begin{array}{rcl}
S_\lambda(\X)\Delta(\X)^{2(k-1)}&=&
\Det\left(S^{\lambda_{i_1}+n-i_1+\dots+n-i_{2k-1}+i_{2k}-n}(\X)\right)_{1\leq i_1,\dots,i_{2k}\leq n}\\
&=&  \Det\left(S^{\lambda_{n-i_1}+i_1+\dots+i_{2k-1}-i_{2k}}(\X)\right)_{0\leq i_1,\dots,i_{2k}\leq n-1}\\
&=&(-1)^{n(n-1)\over 2}\Det\left(S^{\lambda_{n-i_1}+1-n+i_1+\dots+i{2k}}(\X)\right)_{0\leq i_1,\dots,i_{2k}\leq n-1}\\
&=&(-1)^{n(n-1)\over 2}{\cal H}_{{\rm reverse}_n(\lambda)-[(n-1)^n]}(-\overline\X)\\
&=&(-1)^{n(n-1)(k-1)\over 2}{\cal T}_{{\rm
reverse}_n(\lambda)-[((k-1)(n-1))^n]}(-\overline\X).
\end{array}
$$

In particular, from Proposition \ref{genmatsumoto}
$$
\begin{array}{rcl}
R_{[n^{p+(k-1)(n-1)}l]}^{(k),n}(-\overline\X)&=&{\cal
T}^{(k)}_{[((k-1)(n-1))^n]+[p^{n-l}(p+1)^l]}(-\overline \X)\\
&=&(-1)^{n(n-1)(k-1)\over 2}
S_{[(p+1)^lp^{n-l}]}(\X)\Delta(\X)^{2(k-1)}.
\end{array}
$$
But, the Jack polynomial $R_\lambda^{(\alpha),n}$ being homogeneous,
one has
$$
R_\lambda^{(\alpha),n}(-\overline\X)=(-1)^{|\lambda|}R_\lambda^{(\alpha),n}(-\X).
$$
The result follows.
 $\Box$

\begin{remark}\rm
\begin{enumerate}
 \item Note that a
special case of Proposition \ref{Pat-X} appeared in \cite{LT1}.
 \item Proposition \ref{Pat-X} can be reformulated as \\

\noindent The polynomials $P^{(k)}_{n^{p+(k-1)(n-1)}l}(-\X)$ and
$P_{(p^n)+(1^l)}^{(k)}(\X)\Delta(\X)^{2(k-1)}$ are proportional.\\

\noindent This kind of identities relying Jack polynomials in $\X$
and in $-\X$ can be deduced from more general ones involving
Macdonald polynomials when $t$ is specialized to a power of $q$.
This  will be investigated in a forthcoming paper.
\end{enumerate}
\end{remark}
As a special case of Proposition \ref{Pat-X}, the even powers of the
Vandermonde determinants $\Delta(\X)^{2k}$ are Jack polynomials on
the alphabet $-\X$.
\begin{corollary}Setting $l=p=0$ in Equality (\ref{SVand}), one obtains
$$\Delta(\X)^{2k}={(-1)^{(kn(n-1)\over2}\over n!}
\left((k+1)n\atop k+1,\dots,k+1\right)P_{n^{(n-1)k}}^{(k+1)}(-\X).$$
\end{corollary}
In the same way, using Corollary \ref{Q1/XSQ}, one finds a
surprising identity relying Jack polynomials in the alphabets $-\X$
and $\X^\vee$.
\begin{proposition}
One has \begin{eqnarray*}
\Omega^+Q_{\lambda}^{(1/k)}(\X^\vee)\Delta(\X)^{2(k-1)}\Lambda^n(\X)^{n-1}=\\
\frac{(-1)^{{n(n-1)(k-1)\over 2}+|\lambda|}}{n!}\left(nk\atop
k,\dots,k\right)(b_{n^{k(n-1)}}^{(k)})^{-1}Q_{n^{k(n-1)}/\lambda'}^{(k)}(-\X).
 \end{eqnarray*}
\end{proposition}
{\bf Proof} The result is a straightforward consequence of Lemma
\ref{lemma_alt} and Corollary \ref{Q1/XSQ}. $\Box$

 \end{document}